\documentclass[10pt]{article}
\usepackage{latexsym,amsmath,amssymb,epsf}
\makeatletter
\@addtoreset{figure}{section}
\def\thefigure{\thesection.\@arabic\c@figure}
\def\fps@figure{h,t}
\@addtoreset{table}{bsection}

\def\thetable{\thesection.\@arabic\c@table}
\def\fps@table{h, t}
\@addtoreset{equation}{section}

\makeatother

\begin{document}

\newtheorem{theorem}{Theorem}[section]
\newtheorem{definition}[theorem]{Definition}
\newtheorem{lemma}[theorem]{Lemma}
\newtheorem{remark}[theorem]{Remark}
\newtheorem{proposition}[theorem]{Proposition}
\newtheorem{corollary}[theorem]{Corollary}
\newtheorem{example}[theorem]{Example}

\newcommand{\bfi}{\bfseries\itshape}

\newsavebox{\savepar}
\newenvironment{boxit}{\begin{lrbox}{\savepar}
\begin{minipage}[b]{15.5cm}}{\end{minipage}\end{lrbox}
\fbox{\usebox{\savepar}}}

\makeatletter
\title{{\bf The symplectic reduced spaces\\of a Poisson action}}
\author{Juan--Pablo Ortega \\
Institut Nonlin\'eaire de Nice\\
Centre National de la Recherche Scientifique\\
1361, route des Lucioles\\
F-06560 Valbonne, France \\
{\footnotesize Juan-Pablo.Ortega@inln.cnrs.fr}}
\date{3 March 2002}
\maketitle
\makeatother

\addcontentsline{toc}{section}{Abstract}
\begin{abstract} 
During the last thirty years, symplectic or
Marsden--Weinstein reduction has been a major tool in the construction
of new symplectic manifolds and in the study of mechanical systems
with symmetry. This procedure has been traditionally associated to
the canonical action of a Lie group on a symplectic manifold, in the
presence of a momentum map. In this note we show that the symplectic
reduction phenomenon has much deeper roots. More specifically, we
will find symplectically reduced spaces purely within the Poisson
category under hypotheses that do not necessarily imply the
existence of a momentum map. On other words, the right category to
obtain symplectically reduced spaces is that of Poisson manifolds
acted canonically upon by a Lie group.
\end{abstract}

\section{Introduction}

Let $(M , \omega) $ be a symplectic manifold and $G$ be a Lie group
that acts freely and properly on $M$. We will assume that this action
is canonical, that is, it preserves the symplectic form and that it
has an equivariant momentum map $\mathbf{J}: M \rightarrow
\mathfrak{g}^\ast$ associated. Marsden and Weinstein~\cite{mwr} 
showed that for any value $\mu \in \mathbf{J}(M)$ with coadjoint
isotropy subgroup $G _\mu$, the quotient $\mathbf{J}^{-1}(\mu)/G
_\mu$ is a smooth symplectic manifold with a symplectic structure
naturally inherited from that in $M$. This procedure can be reproduced
when, instead of a $\mathfrak{g}^\ast$--valued momentum map, we have a
$G$--valued momentum map in the sense of Alekseev {\it et
al.}~\cite{mcduff 1988, lie group valued maps}.

The study of symplectic reduction in the absence of the freeness
hypothesis on the $G$--action has given rise to the so called {\bfi
Singular Reduction Theory} which has been spelled out over the years
in a series of works. See~\cite{acg, sl, bl, thesis,
cushman sniatycki, hsr}, and references therein.

The first effort to perform symplectic reduction without momentum
maps was carried out in~\cite{optimal} by using the so
called {\bfi optimal momentum   map}. Nevertheless, in the
requirements of the reduction theorem formulated in that paper there
is a ``closedness hypothesis" that is  reminiscent at some level of
the existence of a standard ($\mathfrak{g}^\ast$ or $G$--valued)
momentum map.

In this note we will formulate a symplectic reduction theorem that
does not require this hypothesis and that at the same time works in
the Poisson category. More specifically, we will show that the
Marsden--Weinstein quotients constructed using the (always
available) optimal momentum map  associated to a canonical Lie
group action on the Poisson manifold $(M, \{ \cdot, \cdot\}) $ are
smooth symplectic manifolds, provided that the group action satisfies
a customary properness hypothesis.

\section{The optimal momentum map and the momentum space}

The optimal momentum map was introduced in~\cite{optimal} as a
general method to find the conservation laws associated to the
symmetries of a Poisson system encoded
in the canonical action of a Lie group. We recall its definition.
Let
$(M, \{ \cdot, \cdot\}) $ be a Poisson manifold and $G$ be a Lie
group
that acts properly on
$M$ by Poisson diffeomorphisms via the left action
$\Phi: G
\times M \rightarrow M $. The group of canonical transformations
associated to this action will be denoted by $A_{G}:=\{ \Phi _g: M
\rightarrow M
\mid
 g \in G\}$ and the
canonical projection of $M$ onto the orbit space  by
$\pi_{A_{G}}: M \rightarrow M/ A_{G}=M/G $.
Let $A_{G}'$ be the distribution on $M$ defined by the relation:
\[
A_{G}'(m):=\{ X _f(m)\mid f \in C^\infty(M)^G\}, \quad \quad \mbox{ 
for all} \quad m \in M.
\]
The symbol $X _f$ denotes the Hamiltonian vector field associated to
the function $f \in C ^{\infty}(M)$. Depending on the
context, the distribution
$A_{G}'$ is called the {\bfi $G$--characteristic distribution} or the
{\bfi polar distribution} defined by $A_{G}$~\cite{dual pairs}.
$A_{G}'$ is a smooth integrable  generalized distribution in the
sense of Stefan and Sussman~\cite{stefan, stefan b, sussman}. The
{\bfi optimal momentum map} ${\cal J}$ is defined as the canonical
projection  onto the leaf space of $A_{G}'$, that is,
\[
{\cal J}:M \longrightarrow M / A_{G}'.
\]
By its very definition, the levels sets of ${\cal J}$ are preserved
by the Hamiltonian flows associated to $G$--invariant Hamiltonian
functions and ${\cal J}$ is {\bfi universal} with respect to this
property, that is, any other map whose level sets are preserved by
$G$--equivariant Hamiltonian dynamics factors necessarily through 
${\cal J}$. By construction, the fibers of
${\cal J}$ are the leaves of an integrable generalized
distribution and thereby {\it initial immersed submanifolds} of
$M$~\cite{dazord 1985}. Recall that we say that $N$ is an initial
submanifold of
$M$ when 
the injection $i:N\rightarrow M$ is a smooth immersion that satisfies
that for any manifold
$Z$, a mapping 
$f:Z \rightarrow N $
is smooth iff $i\circ f:Z \rightarrow M $ is smooth. We summarize
this and other elementary properties of the fibers of ${\cal J}$ in
the following proposition.

\begin{proposition}
Let $(M, \{ \cdot, \cdot\}) $ be a Poisson manifold and $G$ be a Lie
group
that acts properly and canonically on
$M$. Let ${\cal J}: M \rightarrow M/ A_{G}'$ be the associated
optimal momentum map. Then for any $\rho\in M / A_{G}'$ we have that:
\begin{description}
\item[ (i)] The level set ${\cal J}^{-1}(\rho)$ is an immersed
initial submanifold of $M$.
\item [(ii)] There is a unique symplectic leaf $\mathcal{L}$ of  $(M,
\{ \cdot, \cdot\}) $ such that $ {\cal J}^{-1}(\rho) \subset
\mathcal{L}$.
\item [(iii)] Let $m \in M$ be an arbitrary element of ${\cal
J}^{-1}(\rho)$. Then, ${\cal J}^{-1}(\rho) \subset M _{G _m}$, with 
$M _{G _m}:=\{ z \in M\mid G _z= G _m\} $.
\end{description}
\end{proposition} 

In the sequel we will  denote by
$\mathcal{L}_\rho$ the unique symplectic leaf of $M$ that contains
${\cal J}^{-1}(\rho)$. Notice that as $\mathcal{L}_\rho$ is also an
immersed initial submanifold of $M $, the injection $i
_{\mathcal{L}_\rho}: {\cal J}^{-1}(\rho) \hookrightarrow
\mathcal{L}_\rho$ is smooth.

The leaf space $M / A_{G}'$ is called the {\bfi  momentum
space} of ${\cal J}$. We will consider it as a topological space
with the quotient topology. Let $m\in M$ be arbitrary such that
${\cal J}(m)=\rho\in M/A_{G}'$. Then, for any $g\in G$, the map
$\Psi_g(\rho)={\cal J}(g\cdot m)\in M/A_{G}'$ defines a continuous
$G$--action on $M/A_{G}'$ with respect to which ${\cal J}$ is
$G$--equivariant. Notice that since this action is not smooth and
$M/A_{G}'$ is not Hausdorff in general, there is no guarantee that
the isotropy subgroups $G _\rho$ are closed, and therefore embedded,
subgroups of $G$. However, there is still something that we can
say:

\begin{proposition}
\label{g rho action}
Let $G _\rho$ be the isotropy subgroup of the element $\rho \in M/
A_{G}'$ associated to the $G$--action on $M/ A_{G}'$ that we just
defined. Then: 
\begin{description}
\item [ (i)] There is a unique smooth structure on $G _\rho$ for
which this subgroup becomes an initial Lie subgroup of $G$ with Lie
algebra $\mathfrak{g}_\rho$ given by 
\[
\mathfrak{g}_\rho=\{ \xi\in\mathfrak{g}\mid  \xi_M(m) \in T _m {\cal
J} ^{-1}(\rho),  \mbox{ for all }  m \in {\cal
J}^{-1}(\rho)\}.
\]
\item [(ii)] With this smooth structure for $G  _\rho$, the left
action
$\Phi^\rho: G _\rho\times
{\cal
J} ^{-1}(\rho)\rightarrow{\cal
J} ^{-1}(\rho)$ defined by
$\Phi^\rho(g, z):= \Phi(g, z)$ is smooth.
\item [(iii)] This action has  fixed isotropies, 
that is, if $z \in {\cal J}^{-1}(\rho)$ then $(G _\rho) _z= G _z$,
and $G _m= G _z$ for all $m \in {\cal J}^{-1}(\rho)$.
\end{description}
\end{proposition}

\noindent\textbf{Proof.\ \ } {\bf (i)} It is a straightforward
corollary of Definition 3 and Proposition 9 in page 290
of~\cite{bourbaki lie 1-3}. Indeed, we can use that result to
conclude the existence of a unique smooth structure for $G _\rho$
with which it becomes an immersed subgroup of $G$ with Lie algebra:
\[
\mathfrak{g}_\rho=\{ \xi \in \mathfrak{g}\mid \mbox{there exists a
smooth curve } c : \Bbb R \rightarrow G _\rho \mbox{ such that } c
(0)=e\mbox{ and } c' (0)= \xi\}.
\] 
An elementary argument shows that
\begin{eqnarray*}
\mathfrak{g}_\rho&=&\{\xi\in \mathfrak{g}\mid  \exp t \xi \cdot m \in
{\cal J}^{-1}(\rho)\mbox{ for all } m \in {\cal J}^{-1}(\rho), t
\in\Bbb R\}\\
 &=&\{ \xi\in\mathfrak{g}\mid  \xi_M(m) \in T _m {\cal
J} ^{-1}(\rho),  \mbox{ for all }  m \in {\cal
J}^{-1}(\rho)\}. 
\end{eqnarray*}

\noindent {\bf (ii)} As ${\cal J}^{-1}(\rho)$ is an initial
submanifold of 
$M$ and $i _\rho \circ \Phi^\rho$ is smooth, with $i _\rho: {\cal
J}^{-1}(\rho)
\hookrightarrow M $ the natural inclusion, then $\Phi^\rho$ is also
smooth. 
 {\bf (iii)} is a straightforward consequence of the
definitions.\quad
$\blacksquare$

\section{The reduction theorem}

We will now introduce our main result. In the statement we will
denote by 
$\pi_\rho: {\cal J}^{-1}(\rho)\rightarrow {\cal J}^{-1}(\rho)/ G
_\rho$ the canonical projection onto the orbit space of the $G
_\rho$--action on ${\cal J}^{-1}(\rho)$ defined in Proposition~\ref{g
rho action}. 

\begin{theorem}[Symplectic reduction by Poisson actions]
\label{Symplectic reduction by Poisson actions}
Let $(M, \{ \cdot, \cdot\}) $ be a smooth Poisson manifold and $G$ be
a Lie group acting canonically and properly on $M$. Let ${\cal J}:M
\rightarrow M/ A_{G}'$ be the optimal momentum map associated to
this action. Then, for any $\rho \in M / A_{G}'$ whose isotropy
subgroup $G _\rho$ acts properly on ${\cal J}^{-1}(\rho)$, the
orbit space $M _\rho:={\cal J}^{-1}(\rho)/ G _\rho$ is a smooth
symplectic regular quotient manifold with symplectic form
$\omega_\rho$ defined by:
\begin{equation}
\label{symplectic 1}
\pi_\rho^\ast\omega_\rho(m)(X _f(m), X _h(m))=\{f, h\} (m), \quad
\mbox{ for any } m \in \mathcal{J}^{-1}(\rho)  \mbox{ and any } f,h
\in C^\infty (M)^G.
\end{equation} 
\end{theorem}

\medskip

\begin{remark}
\normalfont
Let $i_{\mathcal{L}_\rho}: {\cal J}^{-1}(\rho) \hookrightarrow
\mathcal{L}_\rho$ be the natural smooth injection of ${\cal
J}^{-1}(\rho)$ into the symplectic leaf $(\mathcal{L}_\rho,
\omega_{\mathcal{L}_\rho})$ of $(M, \{ \cdot, \cdot\}) $ in which it
is sitting. As $\mathcal{L}_{\rho}$ is an initial submanifold of $M$,
the injection $i _{\mathcal{L}_{\rho}} $ is a smooth map. The form
$\omega_\rho$ can also be written in terms of the symplectic
structure of the leaf
$\mathcal{L}_\rho$  as
\begin{equation}
\label{symplectic 2}
\pi_\rho^\ast\omega_\rho=i _{\mathcal{L}_\rho}^\ast
\omega_{\mathcal{L}_\rho}. 
\end{equation}
In view of this remark we can obtain the standard
Symplectic Stratification Theorem of Poisson manifolds as a
straightforward corollary of Theorem~\ref{Symplectic reduction by
Poisson actions} by taking the group $G=\{ e\}$. In that case the
distribution $ A_{G}'$ coincides with the characteristic
distribution of the Poisson manifold and the level sets of the
optimal momentum map, and thereby the symplectic quotients $M _\rho$,
are exactly the symplectic leaves.  We explicitly point
this out in our next statement.
\quad $\blacklozenge$ 
\end{remark}

\begin{corollary}[Symplectic Stratification Theorem]
Let $(M, \{ \cdot, \cdot\}) $ be a smooth Poisson manifold. Then, $M$
is the disjoint union of the maximal integral leaves of the
integrable distribution $D$ given by
\[
D(m):=\{ X _f(m) \mid f \in C^\infty(M)\}, \quad m \in M.
\]
These leaves are symplectic initial submanifolds  of $M$.
\end{corollary}

\begin{remark}
\normalfont
The only extra hypothesis in the statement of Theorem~\ref{Symplectic
reduction by Poisson actions} with respect to the hypotheses used
in the classical reduction theorems is the properness of the $G
_\rho$--action on ${\cal J}^{-1}(\rho)$. The next example will
show that this is a real hypothesis in the sense that the
properness of the $G
_\rho$--action is not automatically inherited from the properness
of the $G$--action on $M$, as it used to be the case in the presence
of a standard momentum map (see~\cite{optimal}). From this
reduction point of view we can think of the presence of a standard
momentum map as an extra integrability feature of the
$G$--characteristic distribution that makes its integrable leaves
imbedded (and not just initial) submanifolds of $M$ and their
isotropy subgroups automatically closed. \quad $\blacklozenge$ 
\end{remark}

\begin{example}
\normalfont
{\bf On the properness of the $G _\rho$--action.}
As we announced in the previous remark, we now present a situation
where the  $G _\rho$--action on ${\cal
J}^{-1}(\rho)$ is not proper while the 
$G $--action on $M$ satisfies this condition. Let
$M:=\Bbb T^2\times
\Bbb T^2$ be the product of two two--tori whose elements we will
denote by the four--tuples
$(e ^{i \theta_1},e ^{i \theta_2},e ^{i \psi_1},e ^{i \psi_2}) $. We
endow $M$ with the symplectic structure $\omega$ defined by
$\omega:= \mathbf{d} \theta_1\wedge \mathbf{d} \theta_2+ \sqrt{2}\,
\mathbf{d} \psi_1\wedge \mathbf{d} \psi_2$. We now consider the
canonical two--torus action given by $(e ^{i \phi_1},e ^{i \phi_2})
\cdot (e ^{i
\theta_1},e ^{i \theta_2},e ^{i \psi_1},e ^{i
\psi_2}):=(e ^{i (\theta_1+ \phi_1)},e ^{i (\theta_2+ \phi_2)},e ^{i
(\psi_1+
\phi_1)},e ^{i
(\psi_2+ \phi_2)})$. First of all, notice that since the two--torus is
compact this action is necessarily proper. Moreover, as $\Bbb T^2$ 
acts freely, the corresponding orbit space
$M /A _{\Bbb T^2}$ is a smooth manifold such that the projection
$\pi_{A _{\Bbb T^2}}:M \rightarrow M /A _{\Bbb T^2} $ is a surjective
submersion. The polar distribution $A _{\Bbb T^2}' $ does not have
that property. Indeed, 
$C ^{\infty}(M)^{\Bbb T^2}$ comprises all the functions
$f$ of the form $f \equiv f (e ^{i
(\theta_1-\psi_1)},e ^{i
(\theta_2-\psi_2)})$. An inspection of the Hamiltonian flows
associated to such functions readily shows that the leaves of  $A
_{\Bbb T^2}' $, that is, the level sets of the optimal momentum
map ${\cal J}$, are the products of two leaves of an irrational
foliation in a two--torus. Moreover, it can be checked that for
any
$\rho \in M/ A _{{\cal T}^2}' $, the isotropy subgroup $\Bbb
T^2_\rho$ is the product of two discreet subgroups of $S ^1$, each of
which fill densely the circle. We can use this density property to
show that the
$\Bbb T _\rho$--action on ${\cal J}^{-1}(\rho)$ is not proper. Let
$\{(e ^{i \tau_n},e ^{i \sigma_n})\}_{n \in \Bbb N}$ be a
strictly monotone sequence of elements in $\Bbb T^2_\rho$ that
converges to
$(e ,e )$ in $ \Bbb T^2$. Then, for any  sequence $\{ z
_n\}_{n \in\Bbb N} \subset {\cal J}^{-1}(\rho)$ such that $z _n
\rightarrow z \in{\cal J}^{-1}(\rho)$ in ${\cal J}^{-1}(\rho)$ we
have that $(e ^{i \tau_n},e ^{i \sigma_n}) \cdot z _n\rightarrow z $
in ${\cal J}^{-1}(\rho)$. However, since $\Bbb T^2 _\rho$ is endowed
with the discrete topology and $\{(e ^{i \tau_n},e ^{i
\sigma_n})\}_{n \in \Bbb N}$ is strictly monotone it has no
convergent subsequences, which implies that $G _\rho$ does not act
properly on ${\cal J}^{-1}(\rho)$.
 \quad
$\blacklozenge$
\end{example}

\begin{example}
\normalfont
A simplified version of the previous example provides a situation
where the hypotheses of Theorem~\ref{Symplectic reduction by Poisson
actions} are satisfied while all the standard reduction theorems fail.
Namely, there are no momentum maps for this action and, moreover, the
``closedness hypothesis" in~\cite{optimal} is not satisfied. 

Let
$M:=\Bbb T^2\times
\Bbb T^2$ with the same symplectic structure that we had in the
previous example.  We now consider the
canonical circle action given by $e ^{i \phi} \cdot (e ^{i
\theta_1},e ^{i \theta_2},e ^{i \psi_1},e ^{i \psi_2}):=(e ^{i
(\theta_1+ \phi)},e ^{i \theta_2},e ^{i (\psi_1+ \phi)},e ^{i
\psi_2})$.  In this case, 
$C ^{\infty}(M)^{S ^1}$ comprises all the functions
$f$ of the form $f \equiv f (e ^{i \theta_2},e ^{i \psi_2}, e ^{i
(\theta_1-\psi_1)})$. An inspection of the Hamiltonian flows
associated to such functions readily shows that the leaves of  $A _{S
^1}' $, that is, the level sets ${\cal J}^{-1}(\rho)$ of the optimal
momentum map ${\cal J}$, are the product of a two--torus with a leaf
of an irrational foliation (Kronecker submanifold) of another
two--torus. Obviously this is not compatible with the existence of a
($\Bbb R^2$ or
$\Bbb T^2$--valued) momentum map or with the closedness hypothesis
in~\cite{optimal}. Nevertheless, the isotropies $S ^1 _\rho$ coincide
with the circle $S ^1$, whose compactness guarantees that its action
on ${\cal J}^{-1}(\rho)$ is proper. Theorem~\ref{Symplectic reduction by Poisson
actions} automatically guarantees that the quotients of the form
\[
M _\rho:= {\cal J}^{-1}(\rho)/ S ^1_\rho\simeq \left(S ^1\times
_{S^1} S ^1\right)
\times \text{\{Kronecker submanifold of $\Bbb T^2$\}}.
\] 
are symplectic.
 \quad
$\blacklozenge$
\end{example}

\begin{example}
\normalfont
{\bf A Poisson example.} We now use Theorem~\ref{Symplectic reduction
by Poisson actions} to carry out the symplectic reduction of a Poisson
symmetric manifold that was already used in~\cite{optimal} to
illustrate the construction of the optimal momentum map. 
Let
$({\Bbb R}^3,
\{\cdot, \cdot\})$ be the Poisson manifold formed by the Euclidean
three dimensional space ${\Bbb R}^3$ together with the Poisson
structure induced by the Poisson tensor $B$ that in Euclidean
coordinates takes the form:
\[
B= \left(
\begin{array}{ccc}
0 &1 &0\\
-1 &0 &1\\
0 &-1 &0
\end{array}\right).
\]
Consider the action of the additive group $(\Bbb R, +)$ on ${\Bbb
R}^3$ given by $\lambda\cdot(x,y,z):=(x+ \lambda,y,z)$, for any
$\lambda\in\Bbb R$ and any $(x,y,z)\in {\Bbb R}^3$.  This action is
proper and, as we saw in~\cite{optimal}  it does not have a
standard associated momentum map. Nevertheless, it is a Poisson
action and it has an optimal momentum map ${\cal J}$ associated to it
given by the expression
\[
\begin{array}{cccc}
{\cal J}:& {\Bbb R}^3 & \longrightarrow&\Bbb R\\
	&(x,y,z )& \longmapsto & x+z. 
\end{array}
\]
The level sets ${\cal J}^{-1}(c)$ of the optimal momentum map are the
planes given by the equation $x+z=c $ and the isotropy subgroups $
\Bbb R _c $ are always trivial. Therefore, Theorem~\ref{Symplectic reduction
by Poisson actions} concludes that the planes of the
form $x+z=c $ are symplectic submanifolds  of 
the Poisson manifold
$({\Bbb R}^3,
\{\cdot, \cdot\})$. Actually, it is easy to verify that these planes
constitute its symplectic leaves. \quad $\blacklozenge$
\end{example}

\medskip

\noindent\textbf{Proof of the theorem.\ \ } Since by hypothesis the $G
_\rho$--action on ${\cal J}^{-1}(\rho)$ is proper and by
Proposition~\ref{g rho action} it has fixed isotropies,  the
quotient
${\cal J}^{-1}(\rho)/ G _\rho$ is therefore a smooth manifold, and the
projection $\pi_\rho: {\cal J}^{-1}(\rho)\rightarrow {\cal
J}^{-1}(\rho)/ G _\rho$ is a smooth surjective submersion.

We start the proof of the symplecticity of $M _\rho$ by showing
that~(\ref{symplectic 1}) is a good definition for the form
$\omega_\rho$ in the quotient $M _\rho$. Let $m, m' \in {\cal
J}^{-1}(\rho)$ be such that
$\pi_\rho (m)=\pi_\rho (m')$, and $v, w \in T _m {\cal
J}^{-1}(\rho)$, $v', w' \in  T _{m'} {\cal
J}^{-1}(\rho)$ be such that $T _m \pi_\rho \cdot v= T _{m'} \pi_\rho
\cdot v'$, $T _m \pi_\rho \cdot w= T _{m'} \pi_\rho
\cdot w'$. Let $f, f', g, g' \in C^\infty(M)^G$ be such
that $v= X _f(m)$, $v'= X _{f'}(m')$, $w= X _g(m)$, $w'= X
_{g'}(m')$. The condition $\pi_\rho (m)=\pi_\rho (m')$ implies the
existence of an element $k \in G _\rho$ such that $m'= \Phi^\rho_k
(m)$. We also have that $T _m \pi_\rho= T _{m'} \pi_\rho \circ T _m
\Phi^\rho_k$. Analogously, because of the equalities $T _m \pi_\rho
\cdot v= T _{m'} \pi_\rho
\cdot v'$, $T _m \pi_\rho \cdot w= T _{m'} \pi_\rho
\cdot w'$ there exist $G$--invariant functions $h ^1, h ^2 \in
C^\infty(M)^G$ and elements $\xi^1, \xi^2 \in \mathfrak{g}_\rho$ such
that 
\begin{eqnarray*}
X _{f'} (m')- T _m \Phi ^\rho_k \cdot X _f (m)&=& \xi^1_{{\cal
J}^{-1}(\rho)}(m')=X _{h ^1}(m'),\\
X _{g'} (m')- T _m
\Phi ^\rho_k \cdot X _g (m)&=& \xi^2_{{\cal J}^{-1}(\rho)}(m')=X _{h
^2}(m'),
\end{eqnarray*}
or, analogously
\[
X _{f'} (m')=X_{h ^1+f \circ \Phi _{k
^{-1}}}(m')=X_{h ^1+f }(m'),\quad\text{and}\quad  X _{g'} (m')=X_{h
^2+g \circ \Phi _{k ^{-1}}}(m')=X_{h
^2+g }(m').
\]
Hence, we can write
\begin{eqnarray*}
\omega_\rho(\pi_\rho(m'))(v',w')&=&\{f', g'\}(m')=\{h ^1+f ,h ^2+g
\}(m')
	=\{h ^1+f ,h ^2+ g \}(m)\\
	&=&\{f,g\} (m)+\{f ,h ^2  \}(m)+\{h ^1 ,g \}(m)+\{h ^1 ,h ^2
 \}(m)\\
	&=&\{f,g\} (m)+ \mathbf{d}f (m) \cdot \xi^2_{{\cal
J}^{-1}(\rho)}(m)-\mathbf{d}(g+h ^2) (m) \cdot
\xi^1_{{\cal J}^{-1}(\rho)}(m)\\
	&=& \{f,g\} (m)=\omega_\rho(\pi_\rho(m))(v,w).
\end{eqnarray*}
Consequently, $\omega_\rho$ is a well defined two--form on the
quotient
$M _\rho$. Given that $\pi_\rho$ is a smooth surjective submersion,
the form $\omega_\rho$ is clearly smooth. The Jacobi identity for
the bracket $\{ \cdot, \cdot\} $ on $M$ implies that $\omega _\rho$
is closed. These two features of the form $\omega_\rho$ can also be
immediately read out of expression~(\ref{symplectic 2}), whose
equivalence with~(\ref{symplectic 1})   is straightforward. 

It only remains to be shown that $\omega_\rho$ is non degenerate. We
start our argument with a few notations and remarks. Let
$H\subset G$ be the isotropy subgroup of all the elements in ${\cal
J}^{-1}(\rho)$ with respect to the smooth $G _\rho$--action on this
manifold. Recall that by Proposition~\ref{g rho action} this
isotropy subgroup coincides with an isotropy  of  the $G$--action on
$M $. Since by hypothesis the $G$--action on $M$ is proper, the
subgroup
$H\subset G _\rho$ is necessarily compact. Moreover, the Slice
Theorem guarantees that for any point $m \in {\cal J}^{-1}(\rho)$,
there is a $G$--invariant neighborhood $U$ of $m$ in $M$ that is
$G$--equivariantly diffeomorphic to the twist product $G \times_H V
_r$, where $V _r $ is a ball of radius $r$ around the origin in some
vector space $V$ on which $H$ acts linearly.

Let $m \in {\cal J}^{-1}(\rho)$. Suppose that the vector $X _f (m) $,
$f
\in C^\infty(M)^G$, is such that  
\begin{equation}
\label{non degeneracy condition}
\pi_\rho^\ast\omega_\rho(m)(X _f(m), X _h(m))=\{f, h\} (m)=0, \quad
\mbox{ for all }  h \in C^\infty (M)^G.
\end{equation}
In order to prove that $\omega_\rho$ is non degenerate we have to show
that  $X _f (m) \in T _m(G _\rho \cdot m)$. We will do so by using
the local coordinates around the point $m$ provided by the Slice
Theorem. First of all, as $f$ is $G$--invariant $X _f (m) \in T _m M
_H$. Hence, as in local coordinates $M _H\simeq N (H) \times_H V
_r^H$, we have that $X _f (m)=T_{(e, 0)} \pi \cdot(\zeta, v)$, where
$\pi:G \times V_r\rightarrow G \times_H V_r$ is the natural
projection, $\zeta \in {\rm Lie}(N (H)) $, and $v \in V  ^H$. We
recall that $V  ^H$ denotes the fixed points in $V $ by the
action of $H$.

We now rephrase in these local coordinates the condition in~(\ref{non
degeneracy condition}). Indeed, the fact that   
\[\pi_\rho^\ast
\omega_\rho(m)(X _f(m), X _h(m))=\{f, h\} (m)=- \mathbf{d} h (m)
\cdot X _f (m)=0,
\]
for all $h \in C^\infty (M)^G$ amounts to saying that $\mathbf{d}g
(0) \cdot v=0 $ for all the functions $g \in C ^{\infty} (V_r)^H$. On
other words, $v \in \left(\{\mathbf{d}g
(0)\mid  g \in C ^{\infty} (V_r)^H\}\right)^{\circ}$. A known fact
about proper group actions (see Proposition 3.1.1 in~\cite{thesis} or
Proposition 2.14 in~\cite{optimal}) implies that $v \in ((V ^\ast)
^H) ^{\circ}$. Consequently, $v \in V  ^H \cap ((V ^\ast)
^H) ^{\circ}$. We now show that this intersection is trivial and
therefore $v = 0 $ necessarily.

We start by recalling (see again the references that we just
quoted) that the restriction to $(V ^\ast)^H$ of the dual
map associated to the inclusion $i _{V ^H}: V ^H \hookrightarrow V $
is a $H$--equivariant isomorphism from $(V ^\ast) ^H$ to $(V ^H)
^\ast$. Now, as $v \in V  ^H \cap ((V ^\ast)
^H) ^{\circ}$ we have that $\langle \alpha, v \rangle _V=0 $ for
every $\alpha \in ((V )^\ast)^H$. The symbol $\langle \cdot, \cdot
\rangle _V$ denotes the natural pairing of $V$ with its dual. We can
rewrite this condition as
\[
0=\langle \alpha, v \rangle _V=\langle \alpha, i _{V ^H}(v) \rangle
_V=\langle i _{V ^H}^\ast(\alpha), v \rangle _{V^H}.
\]
As the restriction $i _{V ^H} ^\ast| _{(V ^\ast)^H}$ is an
isomorphism, the previous identity is equivalent to $\langle \beta, v
\rangle _{V^H}=0$ for all $\beta \in (V ^H) ^\ast$. Consequently, $v=0
$, as required.

We conclude our argument by noting that as $X _f (m)=T_{(e, 0)} \pi
\cdot(\zeta, 0)$, we have that $X _f (m) \in T _m(G  \cdot m)
\cap A_{G}' (m)= T _m(G _\rho \cdot m) $, which proves the
non degeneracy of  $\omega_\rho$. \quad $\blacksquare$

\medskip

\noindent\textbf{Acknowledgments} This
research was partially supported by the European Commission through
funding for the Research Training Network
\emph{Mechanics and Symmetry in Europe} (MASIE).

\end{document}